\documentclass[11pt]{article}
\usepackage{psfig}
\usepackage{amssymb}
\newcommand{\be}{\begin{equation}}
\newcommand{\ee}{\end{equation}}
\newcommand{\bea}{\begin{eqnarray}}
\newcommand{\eea}{\end{eqnarray}}
\newcommand{\binom}[2]{{#1 \choose #2}}

\begin{document}


\title{An Explicit Formula for Restricted Partition Function through
Bernoulli Polynomials}
\author{Boris Y. Rubinstein
\\
Department of Mathematics, University of California, Davis,
\\One Shields Dr., Davis, CA 95616, U.S.A.}
\date{\today}

\maketitle
\begin{abstract}
Explicit expressions for restricted partition
function $W(s,{\bf d}^m)$ and its quasiperiodic components
$W_j(s,{\bf d}^m)$ (called {\em Sylvester waves}) for a set of positive
integers ${\bf d}^m = \{d_1, d_2, \ldots, d_m\}$ are derived.
The formulas are represented in a form of a finite sum over
Bernoulli polynomials of higher order with periodic coefficients.
\end{abstract}
\section{Introduction}
\label{intro}
The problem of partitions of positive integers has long history started from
the work of Euler who laid a foundation of the theory of partitions 
\cite{GAndrews}, introducing the idea of generating functions.
Many prominent mathematicians contributed to the development of the theory 
using this approach.

J.J. Sylvester provided a new insight and made a
remarkable progress in this field. He found \cite{Sylv1,Sylv2} the 
procedure for computation of a {\it restricted} partition functions, 
and described symmetry properties of such functions. The restricted
partition function $W(s,{\bf d}^m) \equiv W(s,\{d_1,d_2,\ldots,d_m\})$ is a
number of partitions of $s$ into positive integers  $\{d_1,d_2,\ldots,d_m\}$,
each not greater than $s$. The generating function for $W(s,{\bf d}^m)$ 
has a form
\be
F(t,{\bf d}^m)=\prod_{i=1}^m\frac{1}{1-t^{d_{i}}}
 =\sum_{s=0}^{\infty} W(s,{\bf d}^m)\;t^s\;,
\label{genfunc}
\ee
where $W(s,{\bf d}^m)$ satisfies the basic recursive relation
\be
W(s,{\bf d}^m) - W(s-d_m,{\bf d}^m) = W(s,{\bf d}^{m-1})\;.
\label{SW_recursion}
\ee
Sylvester also proved the statement about splitting of the partition 
function into periodic and non-periodic parts and showed that the 
restricted partition function may be presented as a sum of "waves", which 
we call the {\em Sylvester waves}
\be
W(s,{\bf d}^m) = \sum_{j=1} W_j(s,{\bf d}^m)\;,
\label{SylvWavesExpand}
\ee
where summation runs over all distinct factors in the set ${\bf d}^m$.
The wave $W_j(s,{\bf d}^m)$ is a quasipolynomial in $s$ 
closely related to prime roots $\rho_j$ of unit.
Namely, Sylvester showed in \cite{Sylv2} that the wave
$W_j(s,{\bf d}^m)$ is a coefficient of
${t}^{-1}$ in the series expansion in ascending powers of $t$ of
\be
F_j(s,t)=\sum_{\rho_j} \frac{\rho_j^{-s} e^{st}}{\prod_{k=1}^{m}
        \left(1-\rho_j^{d_k} e^{-d_k t}\right)}\;.
\label{generatorWj}
\ee
The summation is made over all prime roots of unit
$\rho_j=\exp(2\pi i n/j)$ for $n$ relatively prime to $j$
(including unity) and smaller than $j$.
This result is
just a recipe for calculation of the partition function and it
does not provide explicit formula.

Using the Sylvester recipe we found in \cite{Rama03} an explicit formula for the
Sylvester wave $W_j(s,{\bf d}^m)$ in a form of finite sum of the Bernoulli
polynomials of higher
order \cite{bat53,NorlundMemo} multiplied by a periodic function of integer 
period $j$. The periodic factor is expressed through the
Eulerian polynomials of higher order \cite{Carlitz1960}.
In this note we show that
it is possible to express the partition function through the Bernoulli
polynomials only.

A special symbolic technique is developed in the theory of polynomials of higher
order, which significantly simplifies computations performed with these
polynomials. A short description of this technique required
for better understanding of this paper is given in the Appendix.

\section{Polynomial part of partition function and \\Bernoulli polynomials}
\label{1}
Consider a polynomial part of the partition function
corresponding to the wave $W_1(s,{\bf d}^m)$. It may be found as a residue of
the generator
\be
F_1(s,t) = \frac{e^{st}}{\prod_{i=1}^m (1-e^{-d_i t})}\;.
\label{generatorW1}
\ee
Recalling the generating function for the Bernoulli polynomials of higher 
order \cite{bat53}:
$$
\frac{e^{st} t^m \prod_{i=1}^m d_i}{\prod_{i=1}^m (e^{d_it}-1)} =
\sum_{n=0}^{\infty} B^{(m)}_n(s|{\bf d}^m)
\frac{t^{n}}{n!}\;,
$$
and a transformation rule
$$
B^{(m)}_n(s|-{\bf d}^m) = B^{(m)}_n(s+\sum_{i=1}^m d_i|{\bf d}^m)\;,
$$
we obtain the relation
\be
\frac{e^{st}}{\prod_{i=1}^m (1-e^{-d_it})} =
\frac{1}{\pi_m} \sum_{n=0}^{\infty} B^{(m)}_n(s+s_m|{\bf d}^m)
\frac{t^{n-m}}{n!}\;,
\label{genfuncBernoulli}
\ee
where
$$
s_m = \sum_{i=1}^m d_i, \ \ \pi_m =  \prod_{i=1}^m d_i\;.
$$

It is immediately seen from (\ref{genfuncBernoulli}) that the coefficient 
of $1/t$ in  (\ref{generatorW1}) is given by the term with $n=m-1$
\be
W_1(s,{\bf d}^m) =
\frac{1}{(m-1)!\;\pi_m}
B_{m-1}^{(m)}(s + s_m | {\bf d}^m)\;.
\label{W_1}
\ee
The polynomial part also admits a symbolic form frequently used in theory of
higher order polynomials
$$
W_1(s,{\bf d}^m) = \frac{1}{(m-1)!\;\pi_m}
\left(s+s_m + \sum_{i=1}^m d_i \;{}^i\! B\right)^{m-1}\;,
$$
where after expansion powers $r_i$ of ${}^i\! B$ are converted into orders
of the Bernoulli numbers
$
{}^i \! B^{r_i} \Rightarrow B_{r_i}
$. 
It is easy to recognize in (\ref{W_1}) the explicit
formula reported recently in
\cite{Beck}, which was obtained by a straightforward computation of the complex 
residue of the generator (\ref{generatorW1}).

Note that basic recursive relation for the Bernoulli polynomials
\cite{NorlundMemo}
$$
B_{n}^{(m)}(s + d_m | {\bf d}^m) -
B_{n}^{(m)}(s | {\bf d}^m) =
n d_m B_{n-1}^{(m-1)}(s | {\bf d}^{m-1})
$$
naturally leads to the basic recursive relation for the polynomial part of
the partition function:
$$
W_1(s,{\bf d}^m) - W_1(s-d_m,{\bf d}^m) =
W_1(s,{\bf d}^{m-1})\;,
$$
which coincides with (\ref{SW_recursion}). This indicates that the 
Bernoulli polynomials of higher order represent a natural
basis for expansion of the partition function and its waves.

\section{Sylvester waves 
and Eulerian polynomials}
\label{j}
Frobenius \cite{Frobenius} studied in great detail the
so-called Eulerian polynomials
$H_n(s,\rho)$ satisfying the generating function
$$
\frac{(1-\rho) e^{st}}{e^t-\rho} = \sum_{n=0}^{\infty} H_n(s,\rho)
\frac{t^n}{n!}, \ \ (\rho \ne 1),
$$
which reduces to definition of the Euler polynomials at fixed value of the
parameter $\rho$
$$
E_n(s) = H_n(s,-1).
$$
The polynomials $H_n(\rho) \equiv H_n(0,\rho)$ satisfy the
symbolic recursion ($H_0(\rho)=1$)
\be
\rho H_n(\rho) = (H(\rho)+1)^n, \ \ \ n>0.
\label{EulerSymbolic}
\ee
The generalization to higher orders is straightforward
$$
\frac{e^{st} \prod_{i=1}^m (1-\rho^{d_i})}
{\prod_{i=1}^m (e^{d_i t}-\rho^{d_i})} =
\sum_{n=0}^{\infty} H_n^{(m)}(s, \rho | {\bf d}^m) \frac{t^n}{n!},
\ \ (\rho^{d_i} \ne 1),
$$
where the corresponding recursive relation for $H_n^{(m)}(s, \rho | {\bf 
d}^m)$ has the form
$$
H_n^{(m)}(s+d_m,\rho | {\bf d}^m)-\rho^{d_m}H_n^{(m)}(s,\rho | {\bf d}^m)=
\left(1-\rho^{d_m}\right)H_n^{(m-1)}(s,\rho | {\bf d}^{m-1})\;.
$$
The Eulerian polynomials of higher order
$H_n^{(m)}(s, \rho | {\bf d}^m)$
introduced by L. Carlitz in \cite{Carlitz1960}
can be defined through the symbolic
formula
$$
H_n^{(m)}(s, \rho | {\bf d}^m) =
\left(s + \sum_{i=1}^{m} d_i \;{}^i\! H(\rho^{d_i}) \right)^n,
$$
where $H_n(\rho)$ computed from the
relation
$$
\frac{1-\rho}{e^t-\rho} = \sum_{n=0}^{\infty} H_n(\rho)
\frac{t^n}{n!},
$$
or using the recursion (\ref{EulerSymbolic}). Using the polynomials 
$H_n^{(m)}(s, \rho | {\bf d}^m)$ we can compute the Sylvester wave of
arbitrary period.

In order to compute the Sylvester wave with period $j>1$ we note that
the summand in the expression (\ref{generatorWj}) can be rewritten as
a product
\be
F_j(s,t) = \sum_{\rho_j}
\frac{e^{st}}{\prod_{i=1}^{k_j} (1-e^{-d_it})} \times
\frac{\rho_j^{-s}}{\prod_{i=k_j+1}^m (1-\rho_j^{d_i}e^{-d_i t})}\;,
\label{generator_product}
\ee
where the elements in ${\bf d}^m$ are sorted in a way that $j$ is a divisor
for first $k_j$ elements (we say that $j$ has weight
$k_j$), and the rest elements in the set are not
divisible by $j$.

Consider the $j$-periodic Sylvester wave $W_j(s,{\bf d}^m)$,
and rewrite the summand in (\ref{generator_product}) as double infinite sum
$$
\frac{\rho_j^{-s}}{\pi_{k_j} \; \prod_{i=k_j+1}^m (1-\rho_j^{d_i})}
\sum_{n=0}^{\infty} B^{(k_j)}_n(s+s_{k_j}|{\bf d}^{k_j})
\frac{t^{n-k_j}}{n!}
\sum_{l=0}^{\infty} H_l^{(m-k_j)}(s_m-s_{k_j}, \rho_j | {\bf
 d}^{m-k_j})
\frac{t^l}{l!},
$$
The coefficient of $1/t$ in the above series is found for
$n+l=k_j-1$, so that we obtain a finite sum:
\bea
W_j(s,{\bf d}^m) & = &
\frac{1}{(k_j-1)! \; \pi_{k_j}}
\sum_{\rho_j}
\frac{\rho_j^{-s}}{\prod_{i=k_j+1}^m (1-\rho_j^{d_i})}
\times \nonumber \\
&&\sum_{n=0}^{k_j-1} \binom{k_j-1}{n}
B^{(k_j)}_n(s+s_{k_j}|{\bf d}^{k_j})
H_{k_j-1-n}^{(m-k_j)}(s_m-s_{k_j}, \rho_j | {\bf
 d}^{m-k_j})\;. 
\label{Wj}
\eea
This expression may be rewritten as a symbolic power
:
\be
W_j(s,{\bf d}^m) =
\frac{1}{(k_j-1)! \; \pi_{k_j}}
\sum_{\rho_j}
\frac{\rho_j^{-s}}{\prod_{i=k_j+1}^m (1-\rho_j^{d_i})}
\left(
s+s_m + \sum_{i=1}^{k_j} d_i \;{}^i\! B + \!\!\!
\sum_{i=k_j+1}^{m} \!\!\! d_i \;{}^i\! H(\rho_j^{d_i})
\right)^{k_j-1} \!\!\!\!\!\!\!\!\;,
\label{Wjsymb}
\ee
which is equal to
\bea
W_j(s,{\bf d}^m) & = &
\frac{1}{(k_j-1)! \; \pi_{k_j}}
\sum_{n=0}^{k_j-1} \binom{k_j-1}{n}
B^{(k_j)}_n(s+s_{m}|{\bf d}^{k_j})
\times \nonumber \\
&& \sum_{\rho_j}
\frac{\rho_j^{-s}}{\prod_{i=k_j+1}^m (1-\rho_j^{d_i})}
H_{k_j-1-n}^{(m-k_j)}[\rho_j |{\bf d}^{m-k_j}]\;,
\label{WjBern}
\eea
where
$$
H_{n}^{(m)}[\rho |{\bf d}^{m}] =
H_{n}^{(m)}(0,\rho |{\bf d}^{m}) =
\left[\sum_{i=1}^{m} d_i \;{}^i\! H(\rho^{d_i})\right]^n\;,
$$
are the Eulerian numbers of higher order
and it is assumed that
$$
H_{0}^{(0)}[\rho | \emptyset] = 1, \
H_{n}^{(0)}[\rho | \emptyset] = 0\;, \ n>0\;.
$$

It should be underlined that the presentation of the Sylvester wave as a
finite sum of the Bernoulli polynomials with periodic coefficients
(\ref{WjBern}) is not unique. The symbolic formula (\ref{Wjsymb}) can be
cast into a sum of the Eulerian polynomials
\bea
W_j(s,{\bf d}^m) & = &
\frac{1}{(k_j-1)! \; \pi_{k_j}}
\sum_{n=0}^{k_j-1} \binom{k_j-1}{n}
B^{(k_j)}_n[{\bf d}^{k_j}]
\times \nonumber \\
&& \sum_{\rho_j}
\frac{\rho_j^{-s}}{\prod_{i=k_j+1}^m (1-\rho_j^{d_i})}
H_{k_j-1-n}^{(m-k_j)}(s+s_m, \rho_j |{\bf d}^{m-k_j})\;,
\label{WjEuler}
\eea
where 
$$
B^{(m)}_n[{\bf d}^{m}] = B^{(m)}_n(0|{\bf d}^{m})
$$
are the Bernoulli numbers of higher order.

It should be noted that the formulae (\ref{WjBern}) and (\ref{WjEuler}) require
summation over all prime roots of unit, and though it is simpler than the
Sylvester recipe using (\ref{generatorWj}), it cannot be considered a completely
explicit formula.

\section{Reduction of Sylvester waves to
Bernoulli polynomials}

A relation between the Eulerian and Bernoulli numbers and polynomials of higher
order was established
in \cite{Carlitz1960}, it may be written as follows:
\bea
&&\frac{(m-1-n)!\;\pi_m\;\rho_j^{s_m}}{(k_j-1-n)!\;\pi_{k_j}\;\rho_j^{s_{k_j}}}
\;\frac{1}{\prod_{i=k_j+1}^m (1-\rho_j^{d_i})}\;
H_{k_j-1-n}^{(m-k_j)}[\rho_j | {\bf d}^{m-k_j}]
\nonumber \\
&=&j^{-(m-k_j)}\left(\prod_{i=k_j+1}^{m} \sum_{r_i=0}^{j-1}
\rho_j^{-d_i r_i}\right)
B_{m-1-n}^{(m-k_j)}
(\!\!\!\sum_{i=k_j+1}^{m} d_i r_i| j{\bf d}^{m-k_j}).
\label{EulerBernrel}
\eea
Using this relation we convert the inner sum in
(\ref{WjEuler}) into
\bea
&&\!\!\!\!\!\!\!\!\!\!\!\!\!\!\!\sum_{\rho_j}
\frac{\rho_j^{-s}}{\prod_{i=k_j+1}^m (1-\rho_j^{d_i})}
H_{k_j-1-n}^{(m-k_j)}(s+s_{m},\rho_j |{\bf d}^{m-k_j})
\nonumber \\
&=&j^{-(m-k_j)}\;
\frac{(k_j-1-n)!\;\pi_{k_j}}{(m-1-n)!\;\pi_m}\;
\sum_{\rho_j}\rho_j^{s_{k_j}-s_m-s}
\nonumber \\
&\times&\left(\prod_{i=k_j+1}^{m} \sum_{r_i=0}^{j-1}
\rho_j^{-d_i r_i}\right)
B_{m-1-n}^{(m-k_j)}
(\!\!\!\sum_{i=k_j+1}^{m} d_i r_i| j{\bf d}^{m-k_j})
\nonumber \\
&=&
j^{-(m-k_j)}\;
\frac{(k_j-1-n)!\;\pi_{k_j}}{(m-1-n)!\;\pi_m}\;
\left(\prod_{i=k_j+1}^{m} \sum_{r_i=0}^{j-1}\right)
\nonumber \\
&\times&
B_{m-1-n}^{(m-k_j)}
(\!\!\!\sum_{i=k_j+1}^{m} d_i r_i| j{\bf d}^{m-k_j})
\bar{\Psi}_j(s+\sum_{i=k_j+1}^{m} d_i (r_i+1)).
\label{innersum}
\eea
Here $\bar{\Psi}_j(s)$
denotes a $j$-periodic {\em prime radical circulator} introduced
in \cite{Cayley1855}
$$
\bar{\Psi}_j(s) = 
\sum_{\rho_j} \rho_j^s.
$$
For prime $j$ it is given by
\be
\bar{\Psi}_j(s) =
	\left\{ \begin{array}{ll}
         \phi(j), & \mbox{$s=0 \pmod{j}$}, \\
	 \mu(j), & \mbox{$s\ne 0 \pmod{j}$},
         \end{array}\right.
\label{primecirc}
\ee
where $\phi(j)$ and $\mu(j)$ denote Euler totient and M\"obius functions.
Considering $j$ as a product of powers of
distinct prime factors
$$
j = \prod_{k} p_{k}^{\alpha_{k}},
$$
one may easily check that for integer values of $s$
\be
\bar{\Psi}_j(s) =
\prod_{k} p_{k}^{\alpha_{k}-1}
\bar{\Psi}_{p_{k}}\left(\frac{s}{p_{k}^{\alpha_{k}-1}}\right),
\label{gencirc}
\ee
where $\bar{\Psi}_k(s) = 0$ for non-integer values of $s$.

It is convenient to introduce a {\it $j$-modified} set of summands
${\bf d}^m_j$
defined as union of subset ${\bf d}^{k_j}$ of
summands divisible by $j$ and
the remaining part ${\bf d}^{m-k_j}$ multiplied by $j$
$$
{\bf d}^m_j = {\bf d}^{k_j} \cup j{\bf d}^{m-k_j},
$$
so that ${\bf d}^m_j$ is divisible by $j$.
Substitution of (\ref{innersum}) into (\ref{WjEuler})
with extension of the outer summation up to $m$
produces the formula for computation of the $j$-periodic
Sylvester wave:
\bea
W_j(s,{\bf d}^m) & = &
\frac{1}{(m-1)! \; \pi_{m}\;j^{m-k_j}}
\left(\prod_{i=k_j+1}^{m} \sum_{r_i=0}^{j-1}\right)
\nonumber \\
&\times&
B^{(m)}_{m-1}(s+s_m+\!\!\!\!\sum_{i=k_j+1}^{m} d_i r_i|
{\bf d}^{m}_j)
\bar{\Psi}_j(s+\!\!\sum_{i=k_j+1}^{m} d_i (r_i+1)).
\label{WjFin}
\eea
The derivation of (\ref{WjFin}) implies that all
terms containing $s$
in powers larger than $k_j-1$
are identically equal to zero.
The polynomial part of the partition function
$W_1(s,{\bf d}^m_j)$ for the $j$-modified set of summands reads:
$$
W_1(s,{\bf d}^m_j) =
\frac{1}{(m-1)! \; \pi_{m}\;j^{m-k_j}}
B^{(m)}_{m-1}(s+s_m+(j-1)\!\!\sum_{i=k_j+1}^{m} d_i|{\bf d}^m_j).
$$
This formula gives rise to the representation of the
$j$-periodic Sylvester wave $W_j(s,{\bf d}^m)$ through the
linear combination of the polynomial part of the
$j$-modified set of summand ${\bf d}^m_j$ multiplied by the
$j$-periodic functions $\bar{\Psi}_j$:
$$
W_j(s,{\bf d}^m)=
\left(\prod_{i=k_j+1}^{m} \sum_{r_i=1}^{j}\right)
W_1(s+\!\!\sum_{i=k_j+1}^{m}\!\! d_i (r_i-j),{\bf d}^m_j)
\bar{\Psi}_j(s+\!\!\sum_{i=k_j+1}^{m}\!\! d_i (r_i-j)),
$$
which is written also as
\be
W_j(s,{\bf d}^m)=
\left(\prod_{i=k_j+1}^{m} \sum_{r_i=0}^{j-1}\right)
W_1(s-\!\!\!\!\sum_{i=k_j+1}^{m}\!\!\!\! d_i r_i,{\bf d}^m_j)
\bar{\Psi}_j(s-\!\!\!\!\sum_{i=k_j+1}^{m}\!\!\!\! d_i r_i).
\label{Wj_as_W1a}
\ee
The last formula shows that each Sylvester wave is
expressed as a linear superposition
of the polynomial parts of
the modified set of summands multiplied by the
corresponding prime circulator. Thus, the formulae
(\ref{primecirc}-\ref{WjFin}) with the
Sylvester splitting formula (\ref{SylvWavesExpand}) provide
a complete explicit solution of the restricted partitions problem.

\section*{Appendix \label{appendix1}}

The symbolic technique for manipulating sums with binomial coefficients by
expanding
polynomials and then replacing powers by subscripts was developed in
nineteenth century by Blissard.
It has been known as symbolic notation and the classical umbral
calculus \cite{Roman1978}.
An example of this notation is also found in \cite{bat53} in
section devoted to the Bernoulli polynomials $B_k(x)$.

The well-known formulas
$$
B_n(x+y) = \sum_{k=0}^{n} {n\choose k} B_k(x) y^{n-k}, \ \
B_n(x) = \sum_{k=0}^{n} {n\choose k} B_k x^{n-k},
$$
are written
symbolically as
$$
B_n(x+y) = (B(x)+y)^n, \ \ B_n(x) = (B+x)^n.
$$
After the expansion the exponents of $B(x)$ and $B$ are converted into the
orders of the Bernoulli polynomial and the Bernoulli number, respectively:
$$
[B(x)]^k \Rightarrow B_k(x), \ \ \ B^k \Rightarrow B_k.
$$
We use this notation in its extended version suggested in
\cite{NorlundMemo} in order to make derivation more clear and
intelligible.
N\"orlund introduced the Bernoulli polynomials of higher order defined 
through
the recursion
$$
B_{n}^{(m)}(x|{\bf d}^m) =
\sum_{k=0}^n  \binom{n}{k} d^k B_k(0) B_{n-k}^{(m-1)}(x|{\bf d}^{m-1}),
$$
starting from $B_{n}^{(1)}(x|d_1) = d_1^n B_n(\frac{x}{d_1})$.
In symbolic notation it takes form
$$
B_{n}^{(m)}(x) =  
\left(
d_m B(0) + B^{(m-1)}(x)
\right)^n,
$$
and recursively reduces to more symmetric expression
$$
B_{n}^{(m)}(x|{\bf d}^m) =
\left(
x + d_1 \;{}^1\! B(0) +
d_2 \;{}^2\! B(0) + \ldots + d_m \;{}^m\! B(0)
\right)^n =
\left(
x + \sum_{i=1}^m d_i \;{}^i\! B(0)
\right)^n,
$$
where each  $[{}^i \! B(0)]^k$ is converted into $B_k(0)$.


\end{document}